\documentclass[a4paper,12pt]{article}
\usepackage[utf8]{inputenc}
\usepackage[american]{babel}
\usepackage{anysize}

\usepackage{amsmath}
\usepackage{amssymb}
\usepackage{latexsym}
\usepackage{mathrsfs}

\usepackage{amsthm}

\usepackage{indentfirst}

\usepackage{color}

\usepackage{ifpdf}

\ifpdf
\usepackage[pdftex]{graphicx}
\else
\usepackage{graphicx}
\fi

\newtheorem{theorem}{Theorem}

\newtheorem{claim}[theorem]{Claim}
\theoremstyle{definition}
\newtheorem{definition}{Definition}

\theoremstyle{remark}
\newtheorem{remark}{Remark}
\newcommand{\conv}{\mathop{\rm conv}\nolimits}
\usepackage[hidelinks]{hyperref}

\begin{document}
\sloppy
\ifpdf
\DeclareGraphicsExtensions{.pdf, .jpg, .tif, .mps}
\else
\DeclareGraphicsExtensions{.eps, .jpg, .mps}
\fi

\title{Linear Duality and Its Nonlinear Extension via Topological Covering Theorems}

\author{Mikhail Bludov }

\date{}

\maketitle

\begin{abstract}

Linear programming duality and core existence theorems in cooperative game theory are closely related. For example, in the 1960s, Bondareva used linear duality to characterize when the core of a game with transferable utility is nonempty. Scarf later established a core existence theorem for balanced games with non-transferable utility. Shapley derived this result from his KKMS covering theorem, demonstrating how topological covering theorems yield core existence results.

In this note we give a new KKM-type covering proof of linear duality. We then extend the argument to inequalities with continuous perturbations. Under some asymptotic growth assumptions on perturbations, we obtain a feasible point and balancing coefficients satisfying complementary slackness, a Lagrangian lower bound, and a Scarf-type feasibility criterion. The latter yields Scarf's theorem for normalized cooperative games with non-transferable utility.  

Using Musin's 2017 KKM-type theorems with boundary conditions and the
author's 2026 results on cover degrees relative to point configurations,
we also obtain a complementary slackness theorem for continuous functions
on compact connected oriented manifolds with nonempty boundary.

\end{abstract}



\section{Introduction}

The linear duality theorem and the Knaster--Kuratowski--Mazurkiewicz (KKM)
lemma are classical results with applications to optimization theory,
discrete geometry, game theory, and topological combinatorics. We recall
both results.

\begin{theorem}[Linear duality]
\label{IntroLinDual}
Let $A$ be a real $m\times n$ matrix, $b\in\mathbb{R}^m$, and
$c\in\mathbb{R}^n$. Consider
\[
 \min_{Ax\geq b}\langle x,c\rangle,\qquad
 \max_{\lambda\geq0,\ A^T\lambda=c}\langle\lambda,b\rangle,
\]
called the primal and dual problems, respectively.
\begin{enumerate}
\item \textit{Weak duality.} For every $x$ with $Ax\geq b$ and every
$\lambda\geq0$ with $A^T\lambda=c$,
\[
 \langle x,c\rangle\geq\langle\lambda,b\rangle.
\]
\item \textit{Strong duality.} If either extremum exists, then both exist
and they are equal.
\end{enumerate}
\end{theorem}

\begin{theorem}[KKM lemma]
\label{IntroKKM}
Let $\Delta=\conv\{v_0,\ldots,v_d\} \subseteq \mathbb{R}^d$ be a simplex, and let
$F_0,\ldots,F_d$ be closed subsets of $\Delta$. Suppose that, for every
nonempty $I\subseteq\{0,\ldots,d\}$,
\[
 \conv\{v_i:i\in I\}\subseteq\bigcup_{i\in I}F_i.
\]
Then $\bigcap_{i=0}^dF_i\neq\emptyset$.
\end{theorem}

The KKM lemma \cite{KKM} is a combinatorial form of Brouwer's fixed point
theorem. Early work on linear duality includes von Neumann's 1947 note
\cite{VonNeumannMaximum}. The relation between linear programming and
matrix games (finite two-player zero-sum games) was studied by Gale,
Kuhn, and Tucker \cite{GaleKuhnTucker} and Dantzig \cite{DantzigGames}.
In particular, linear duality can be derived by reducing linear programs
to matrix games and applying von Neumann's minimax theorem for these
games; see \cite{Adler,VonStengel}
for complete proofs and a discussion of the earlier reductions.
Minimax theorems also admit proofs using fixed point theorems
\cite{Kakutani} and KKM-type covering theorems;
see, for example, \cite{Sion}.
The route through minimax gives an indirect KKM-type proof of linear
duality. In this paper we propose a direct KKM-type argument that yields
several new extensions of the linear duality theorem.

The KKM lemma has many extensions. Komiya
\cite{Komiya} generalized the Knaster--Kuratowski--Mazurkiewicz--Shapley
(KKMS) theorem to polytopes. Musin
\cite{MusH,MusBC} explored connections between covering theorems and
homotopy invariants. In \cite{BludovBalanced}, the author studied the
connection between these invariants and convex representations of
target points by given point configurations.
In \cite{BludovMusin}, the authors applied these methods to
cooperative game theory, introducing a framework for generalized
cooperative games and proving Scarf-type results.

There is a natural connection between linear optimization and cooperative
game theory. The Bondareva--Shapley criterion for nonemptiness of the core
of a cooperative game with transferable utility follows from linear
programming duality \cite{Bon,Sh67}.
Scarf's theorem gives a sufficient condition for nonemptiness of the core
in a more general class of games with non-transferable utility
\cite{Scarf}.
Scarf's theorem can itself be derived from covering theorems such as
Shapley's KKMS theorem \cite{Shapley73}. Following this covering approach,
we obtain a theorem that unifies and generalizes both linear duality
and Scarf's theorem.
We first illustrate the method by proving linear duality.

We begin with some terminology. Fix a real $m\times n$ matrix $A$ and
vectors $b\in\mathbb{R}^m$ and $c\in\mathbb{R}^n$, and consider the
system $Ax\geq b$. Write $a_i$ for the $i$th row of $A$ and
$[m]=\{1,\ldots,m\}$.

\begin{definition}[Feasibility and active constraints]
\label{def:FeasibleActive}
A point $x\in\mathbb{R}^n$ is \emph{feasible} for the system $Ax\geq b$
if it satisfies every inequality. An inequality is \emph{active} at a
feasible point when equality holds. The active index set is
\[
 I(x)=\{i\in[m]:\langle a_i,x\rangle=b_i\}.
\]
\end{definition}

\begin{definition}[Balancing coefficients and complementary slackness]
\label{def:BalanceComplementarity}
Nonnegative coefficients $\lambda_i$, $i\in[m]$, are
\emph{balancing coefficients for $(A,c)$} if
\[
 \sum_i\lambda_i a_i=c.
\]
Their \emph{support} is $\operatorname{supp}\lambda=\{i:\lambda_i>0\}$.
A feasible point $x$ and balancing coefficients $\lambda$ for $(A,c)$ satisfy
\emph{complementary slackness} if
$\operatorname{supp}\lambda\subseteq I(x)$.
\end{definition}

When both constraint systems are feasible, the proof of
Theorem~\ref{IntroLinDual} consists of two ingredients.
First, we need to find a feasible point $x^*$ and balancing coefficients
$\lambda$ for $(A,c)$ satisfying complementary slackness. For any such pair, we have
\[
 \langle x^*,c\rangle
 =\sum_i\lambda_i\langle a_i,x^*\rangle
 =\sum_i\lambda_i b_i.
\]

The second ingredient is weak duality, which shows that $x^*$ attains
the primal minimum and $\lambda$ attains the dual maximum.
The nontrivial step is to find such a pair; we do this using a KKM-type
covering theorem.

We now replace the halfspaces
$\{x:\langle a_i,x\rangle\leq b_i\}$ by sets with nonlinear perturbed boundaries.

Let $A=(a_1,\ldots,a_m)$ be an indexed
family of vectors spanning $\mathbb{R}^n$, and let $h$ be a unit vector with
$\langle a_i,h\rangle>0$ for every $i$. Let $H$ be the hyperplane
orthogonal to $h$, and let $\pi_H$ be the orthogonal projection onto $H$:
\[
 H=h^\perp,\qquad \pi_Hx=x-\langle x,h\rangle h.
\]
Choose $b_i\in\mathbb{R}$ and continuous perturbations
$r_i:H\to\mathbb{R}$, and define
\[
 U_i=\{x\in\mathbb{R}^n:\langle a_i,x\rangle
                  \leq b_i+r_i(\pi_Hx)\}.
\]
Write $\operatorname{cone}(A)=\{\sum_i\lambda_i a_i:\lambda_i\geq0\}$.

Continuity of $r_i$ and positivity of $\langle a_i,h\rangle$ give
\[
 \begin{aligned}
 x\notin\operatorname{int}(U_i)
 &\quad\Longleftrightarrow\quad
 \langle a_i,x\rangle\geq b_i+r_i(\pi_Hx),\\
 x\in\partial U_i
 &\quad\Longleftrightarrow\quad
 \langle a_i,x\rangle=b_i+r_i(\pi_Hx).
 \end{aligned}
\]
Thus a point $x$ is feasible for the perturbed system exactly when
$x\notin\bigcup_i\operatorname{int}(U_i)$.
The $i$th inequality is active at such a point exactly when $x\in\partial U_i$;
complementary slackness requires this whenever $\lambda_i>0$.
No convexity of the sets $U_i$ is assumed.

For a target $c\in\operatorname{cone}(A)$ and a continuous function
$r:H\to\mathbb{R}$, consider the objective
$\langle x,c\rangle-r(\pi_Hx)$ and its sublevel sets
\[
 D_b=\{x\in\mathbb{R}^n:\langle x,c\rangle\leq r(\pi_Hx)+b\}
 \qquad(b\in\mathbb{R}).
\]
No growth condition is imposed on $r$.
In the next theorem, an intersection indexed by an empty support is
understood as $\mathbb{R}^n$.

\begin{theorem}[Linear duality with perturbations]
\label{IntroNonlinear}
Let $c\in\operatorname{cone}(A)$. Assume that
\[
 \max\{r_i(y),0\}=o(\|y\|)\qquad(i\in[m]).
\]
Suppose also that there is a set $J\subseteq\{1,\ldots,m\}$ such that
\[
 r_j(y)=o(\|y\|)\quad(j\in J),\qquad
 \operatorname{cone}\{a_j:j\in J\}=\operatorname{cone}(A).
\]
The asymptotic conditions are taken as $\|y\|\to\infty$ in $H$
and are omitted when $H=\{0\}$.
Then the following assertions hold.
\begin{enumerate}
\item There are
$x^*\in\mathbb{R}^n$ and coefficients $\lambda_i\geq0$, $i\in[m]$, such that
\[
 \begin{gathered}
 x^*\notin\bigcup_{i=1}^m\operatorname{int}(U_i),\qquad
 \sum_{i=1}^m\lambda_i a_i=c,\\
 \lambda_i>0\Longrightarrow x^*\in\partial U_i\qquad(i\in[m]).
 \end{gathered}
\]
The coefficients $(\lambda_i)_{i=1}^m$ can be chosen so that at most $n$ of them are positive.

\item \emph{Lagrangian lower bound.} We have
\[
 \begin{aligned}
 &\inf_{x\notin\bigcup_{i=1}^m\operatorname{int}(U_i)}
       \bigl(\langle x,c\rangle-r(\pi_Hx)\bigr)\\
 &\qquad\geq
 \sup_{\substack{\lambda_i\geq0\ (i\in[m])\\\sum_{i=1}^m\lambda_i a_i=c}}
 \left(\sum_{i=1}^m\lambda_i b_i+
       \inf_{y\in H}\left(\sum_{i=1}^m\lambda_i r_i(y)-r(y)\right)\right).
 \end{aligned}
\]
All infima and suprema here are understood in the extended real sense.
If the left-hand side equals $-\infty$, every inner infimum, and hence
the right-hand side, also equals $-\infty$.

\item \emph{Scarf-type feasibility.} Let $b\in\mathbb{R}$ and suppose
that, for every choice of coefficients $\lambda_i\geq0$, $i\in[m]$, with
$\sum_{i=1}^m\lambda_i a_i=c$,
\[
 \bigcap_{\substack{i\in[m]\\\lambda_i>0}}U_i\subseteq D_b.
\]
Then
\[
 D_b\setminus\bigcup_{i=1}^m\operatorname{int}(U_i)\neq\emptyset.
\]
It suffices to check the inclusion for balancing coefficients for $(A,c)$ with
at most $n$ positive entries.
\end{enumerate}
\end{theorem}

The covering extension theorems of \cite{BludovBalanced,MusH} also apply
to manifolds. We use them to obtain a manifold analogue of the first
assertion of Theorem~\ref{IntroNonlinear}.

Throughout the paper, manifolds denoted by $M$ are assumed to be connected.
Let $M$ be a compact oriented $n$-manifold, $n\geq1$, with nonempty
boundary $N=\partial M$, carrying the induced orientation.
Let $\phi_i:M\to\mathbb{R}$, $i\in[m]$, be continuous and put
\[
 \begin{gathered}
 U_i=\{(x,t)\in M\times\mathbb{R}:t\leq\phi_i(x)\},\\
 \Phi(x)=\max_{i\in[m]}\phi_i(x),\qquad
 I(x)=\{i\in[m]:\phi_i(x)=\Phi(x)\}.
 \end{gathered}
\]
Thus $I(x)$ is the set of active indices at $(x,\Phi(x))$.
A point $(x,t)\in M\times\mathbb{R}$ is feasible if $t\geq\Phi(x)$,
equivalently if it lies
outside every $\operatorname{int}U_i$, with interiors taken relative
to $M\times\mathbb{R}$.

\begin{theorem}[Complementary slackness on a manifold]
\label{ManifoldComplementarity}
Let $v_1,\ldots,v_m\in\mathbb{R}^n$. Suppose there is a continuous map
$g:N\to S^{n-1}$ with $\deg(g)\neq0$ such that
$\langle v_i,g(x)\rangle<0$ for every $x\in N$ and $i\in I(x)$.

Then there are a feasible point $(x^*,t^*)$, with $x^*\in M\setminus N$
and $t^*=\Phi(x^*)$, and coefficients $\lambda_i\geq0$, $i\in[m]$,
such that $\sum_{i=1}^m\lambda_i=1$,
$\sum_{i=1}^m\lambda_i v_i=0$ and
$t^*=\phi_i(x^*)$ whenever $\lambda_i>0$.
\end{theorem}

Bludov and Musin \cite[Theorem~2.6]{BludovMusin} prove a
Borsuk--Ulam-type KKMS theorem and apply it to cooperative games with
threats. The antipodal covering result in
Theorem~\ref{AntipodalCover} yields the following symmetric variant
of Theorem~\ref{IntroNonlinear}.

Let $h\in\mathbb{R}^n$ be a unit vector, $n\geq2$, and put $H=h^\perp$.
Let $p_1,\ldots,p_{2m}$ be the distinct vertices of a centrally symmetric
polytope $P\subset H$, full-dimensional in $H$, with
$p_{i+m}=-p_i$ for $i\in[2m]$, here it is assumed that $i+m$ is taken modulo $2m$. Assume that $P$ is simplicial, meaning
that every proper face is a simplex. Choose $\alpha_i>0$, $i\in[2m]$,
with $\alpha_{i+m}=\alpha_i$, and set
\[
 a_i=\alpha_i(h+p_i)\qquad(i\in[2m]).
\]
Thus $a_i$ and $a_{i+m}$ are reflections of one another across
$\mathbb{R}h$.

Let $f_i^+,f_i^-:H\to\mathbb{R}$ be continuous functions satisfying
$f_i^\pm(y)=o(\|y\|)$ as $\|y\|\to\infty$ for $i\in[2m]$.
Constant terms are included in these functions. Define
\[
 \begin{aligned}
 \phi_i(y)&=\max\{f_i^+(y)-\langle p_i,y\rangle,\,
                       f_i^-(y)+\langle p_i,y\rangle\},\\
 U_i&=\{y+th:y\in H,\ t\leq\phi_i(y)\}
 \qquad(i\in[2m]).
 \end{aligned}
\]
Write $\Phi(y)=\max_{i\in[2m]}\phi_i(y)$. A point $y+th$ is feasible
precisely when $t\geq\Phi(y)$.

\begin{theorem}[Antipodal complementary slackness]
\label{SymmetricPerturbations}

Suppose that, for every $i\in[m]$, the \emph{complementary dominance
condition} holds, meaning that either
\[
 f_i^-\leq f_{i+m}^+,\qquad f_{i+m}^-\leq f_i^+
 \quad\text{on }H,
\]
or
\[
 f_i^+\leq f_{i+m}^-,\qquad f_{i+m}^+\leq f_i^-
 \quad\text{on }H.
\]
For $c= h$, there are $x^*\in\mathbb{R}^n$
and coefficients $\lambda_i\geq0$, $i\in[2m]$, such that
\[
 \begin{gathered}
 x^*\notin\bigcup_{i=1}^{2m}\operatorname{int}(U_i),\qquad
 \sum_{i=1}^{2m}\lambda_i a_i=c,\\
 \lambda_i>0\Longrightarrow x^*\in\partial U_i\qquad(i\in[2m]).
 \end{gathered}
\]
At most $n$ coefficients among $(\lambda_i)_{i=1}^{2m}$ are positive. 
\end{theorem}

The paper is organized as follows. Section~2 presents the covering
theorems used in the proofs. Sections~3 and~\ref{sec:ProofNonlinear}
prove linear duality and its nonlinear extension, Theorem~\ref{IntroNonlinear},
respectively.
Section~\ref{sec:ProofManifold} proves the manifold and symmetric versions,
Theorems~\ref{ManifoldComplementarity} and~\ref{SymmetricPerturbations},
and explains how the former yields the
first assertion of Theorem~\ref{IntroNonlinear} and hence
Theorem~\ref{IntroLinDual}.
Section~\ref{sec:NTUGames} applies Theorem~\ref{IntroNonlinear} to prove
Scarf's theorem for normalized cooperative games with non-transferable
utility.
Section~\ref{sec:Conclusion} summarizes the results.

\section{Topological preliminaries}

Vector inequalities are understood coordinatewise. We use the Euclidean
scalar product $\langle\cdot,\cdot\rangle$ and norm $\|\cdot\|$.
For a finite family of vectors $X$, $\conv(X)$ denotes its convex hull,
$\operatorname{cone}(X)$ its nonnegative span, and $\operatorname{span}(X)$ its
linear span. For a matrix $A$, the convex hull $\conv(A)$ and cone $\operatorname{cone}(A)$ mean those of
its row vectors. Write
$B^d=\{x\in\mathbb{R}^d:\|x\|\leq1\}$, and $S^{d-1}=\partial B^d$.
For a subset $S$ of some Euclidean space, the symbols $\operatorname{int}(S)$, $\partial S$, and
$\overline S$ denote its interior, boundary, and closure in the ambient
Euclidean space.
For $i\in[2m]$, the index $i+m$ is always interpreted modulo $2m$,
with its representative in $[2m]$; in particular, $(i+m)+m=i$.

Let $M$ be a compact connected oriented $d$-manifold with boundary $N$,
and give $S^{d-1}=\partial B^d$ its standard boundary orientation.
For $d\geq2$, the degree of $g:N\to S^{d-1}$ means the total degree
$\deg(g)=\sum_\alpha\deg(g|_{N_\alpha})$, where $N_\alpha$ are the
connected components of $N$, with their induced boundary orientations.
For $d=1$, with $S^0=\{-1,+1\}$, we define
$\deg(g)=\sum_{x\in g^{-1}(+1)}\varepsilon_x$, where
$\varepsilon_x\in\{-1,+1\}$ is the induced orientation sign of the
boundary point $x$.

\begin{definition}[Balanced sets]
\label{def:BalancedIntersection}
Let $V=(v_1,\ldots,v_m)$ be an indexed family of vectors in
$\mathbb{R}^d$, not necessarily distinct, and let $c\in\mathbb{R}^d$.
For a nonempty set $S\subseteq[m]$, both $S$ and the corresponding
subfamily $(v_i)_{i\in S}$ are called \emph{$(V,c)$-balanced} if
$c\in\conv\{v_i:i\in S\}$, equivalently, if there are coefficients
$\mu_i\geq0$, $i\in S$, with
\[
 \sum_{i\in S}\mu_i=1,\qquad \sum_{i\in S}\mu_i v_i=c.
\]
For sets $F_1,\ldots,F_m$ labeled by $V$, a nonempty intersection
$\bigcap_{i\in S}F_i$ is called \emph{$(V,c)$-balanced} when $S$ is
$(V,c)$-balanced. Any point in such an intersection is called a
\emph{$(V,c)$-balanced intersection point}.
\end{definition}

\begin{definition}[Partition of unity]
\label{def:PartitionUnity}
Let $O_1,\ldots,O_m$ be an open cover of a space $X$.
A \emph{partition of unity subordinate to this cover} is a family of
continuous functions $\psi_i:X\to[0,1]$ such that
\[
 \sum_i\psi_i(x)=1\quad(x\in X),\qquad \operatorname{supp}\psi_i\subseteq O_i.
\]
Here $\operatorname{supp}\psi_i$ is the closure in $X$ of
$\{x:\psi_i(x)>0\}$. In particular, $\psi_i(x)>0$ implies $x\in O_i$.
\end{definition}

Every finite open cover of a compact metric space admits a subordinate
partition of unity.

Let $M$ be a compact oriented $d$-manifold, $d\geq1$, with boundary
$N$ carrying the induced orientation. Let $\mathcal C=(C_i)_{i=1}^m$
be a finite closed cover of $N$, with labels
$A=(a_1,\ldots,a_m)$, with $a_i\in\mathbb{R}^d$, and target
$c\in\mathbb{R}^d$. Suppose that $\mathcal C$ has no $(A,c)$-balanced
intersection.

\begin{definition}[Degree of a cover {\cite[Definition~4.3]{BludovBalanced}}]
\label{def:CoverDegree}
Choose open neighborhoods $U_i\supseteq C_i$ in $N$ such that, for every
nonempty $S\subseteq[m]$,
\[
 \bigcap_{i\in S}U_i\neq\emptyset
 \quad\Longleftrightarrow\quad
 \bigcap_{i\in S}C_i\neq\emptyset.
\]
Let $(\psi_i)$ be a partition of unity subordinate to $(U_i)$. Put
\[
 q(x)=\frac{\sum_i\psi_i(x)(a_i-c)}
              {\|\sum_i\psi_i(x)(a_i-c)\|}\qquad(x\in N).
\]
The \emph{degree of the cover relative to $(A,c)$} is
$\deg_{A,c}(\mathcal C)=\deg(q)$.
\end{definition}

This degree is independent of the neighborhoods and the partition
of unity \cite[Lemmas~2.4 and~4.2]{BludovBalanced}. The next covering-extension theorem was proved in~\cite{BludovBalanced}.

\begin{theorem}[Covering extension {\cite[Theorem~3]{BludovBalanced}}]
\label{CoverExtensionManifolds}
Let $M$ be a compact connected oriented $d$-manifold, $d\geq1$, with nonempty
boundary $N$, carrying the induced orientation. Let
$A=(a_1,\ldots,a_m)$, with $a_i\in\mathbb{R}^d$, satisfy $\dim\conv(A)=d$,
and let $c\in\mathbb{R}^d$.
A finite closed cover $\mathcal C=(C_i)_{i=1}^m$ of $N$ with no
$(A,c)$-balanced intersection extends to a closed cover
$(G_i)_{i=1}^m$ of $M$ with $G_i\cap N=C_i$ and
without $(A,c)$-balanced intersections, if and only if
$\deg_{A,c}(\mathcal C)=0$.
\end{theorem}

Theorem~\ref{CoverExtensionManifolds} yields the following Komiya-type
theorems for balls with spherical Voronoi-type boundary conditions.

\begin{theorem}[Komiya-type theorem for a ball]
\label{DualPolKKM}
Let $d\geq1$, and let $A=(a_1,\ldots,a_m)$ be an indexed family of
vectors in $\mathbb{R}^d$ whose convex hull has nonempty interior.
Let $F_1,\ldots,F_m$ be closed subsets of $B^d$ covering $B^d$.
Suppose that, for $u\in S^{d-1}$,
\[
 u\in F_i\quad\Longrightarrow\quad
 \langle u,a_i\rangle=\min_{j\in[m]}\langle u,a_j\rangle.
\]
Then, for every $c\in\conv(A)$, the cover has an $(A,c)$-balanced intersection
point.
\end{theorem}

\begin{proof}
First suppose that $c\in\operatorname{int}\conv(A)$.

Fix $u\in S^{d-1}$ and a nonempty set $S\subseteq[m]$ such that
$u\in\bigcap_{i\in S}F_i$. The boundary condition implies that all
labels $a_i$, $i\in S$, lie in the exposed face
\[
 \left\{x\in\conv(A):\langle u,x\rangle
       =\min_{z\in\conv(A)}\langle u,z\rangle\right\}.
\]
Since $u\neq0$ and $\conv(A)$ has nonempty interior, this face is proper
and lies in $\partial\conv(A)$. Thus
$c\notin\conv\{a_i:i\in S\}$, so $S$ is not $(A,c)$-balanced.
Hence the boundary cover has no $(A,c)$-balanced intersections. Moreover,
because $c\in\operatorname{int}\conv(A)$,
\[
 \langle u,a_i-c\rangle
 =\min_{z\in\conv(A)}\langle u,z\rangle-\langle u,c\rangle<0
 \qquad(i\in S).
\]

Choose a partition of unity $(\psi_i)$ subordinate to a sufficiently
small open enlargement of this boundary cover as in
Definition~\ref{def:CoverDegree}, preserving these strict inequalities,
and put
\[
 g(u)=\sum_i\psi_i(u)(a_i-c)\qquad(u\in S^{d-1}).
\]

Since $\langle u,g(u)\rangle<0$, the normalized straight-line homotopy
joins $g/\|g\|$ to the antipodal map $u\mapsto-u$, which is not
of degree zero. Thus the boundary cover has nonzero degree relative
to $(A,c)$. Theorem~\ref{CoverExtensionManifolds}, applied with
$M=B^d$, labels $A$, and target $c$, therefore implies that its extension
$(F_i)_{i=1}^m$ to $B^d$ has an $(A,c)$-balanced intersection point.

Finally, when $c\in\partial\conv(A)$, choose
$c_k\in\operatorname{int}\conv(A)$ with $c_k\to c$.
The interior case gives nonempty index sets $S_k$ with
$c_k\in\conv\{a_i:i\in S_k\}$ and
$\bigcap_{i\in S_k}F_i\neq\emptyset$.
Passing to a subsequence, we may assume that $S_k=S$ is fixed.
Its intersection is nonempty, and the closed convex hull of its
labels contains the limit $c$.
\end{proof}

For a fixed target $c$, the same proof gives the conclusion under
the \emph{strict boundary condition}
\[
 u\in F_i\cap S^{d-1}
 \quad\Longrightarrow\quad
 \langle u,a_i-c\rangle<0,
\]
without requiring the minimum condition in the theorem.
In this case every $(A,c)$-balanced intersection point lies in the open ball.

The next theorem gives an antipodal analogue of
Theorem~\ref{DualPolKKM}, in the spirit of
\cite[Theorem~2.6]{BludovMusin}.

\begin{theorem}[Antipodally relabeled Komiya-type theorem for a ball]
\label{AntipodalCover}
Let $V=(v_1,\ldots,v_{2m})$ be the indexed family of distinct vertices of a full-dimensional
centrally symmetric simplicial polytope $P\subset\mathbb{R}^d$,
$d\geq1$, with $v_{i+m}=-v_i$ for $i\in[2m]$.
Let $(G_i)_{i=1}^{2m}$ be a closed cover of $B^d$ such that
\[
 u\in G_i\cap S^{d-1}\quad\Longrightarrow\quad
 \langle v_i,u\rangle=\min_{j\in[2m]}\langle v_j,u\rangle.
\]

Let $T:[2m]\to[2m]$ be a permutation such that, for every $i\in[m]$,
either $T(i)=i$ and $T(i+m)=i+m$, or
$T(i)=i+m$ and $T(i+m)=i$. Thus $T$ fixes or exchanges the two indices in each pair and satisfies
$T\circ T=\operatorname{id}$.

Then there exists a nonempty $(V,0)$-balanced set $S\subseteq[2m]$ such that
\[
 \operatorname{int}(B^d)\cap\bigcap_{i\in S}G_{T(i)}\neq\varnothing.
\]
\end{theorem}

\begin{proof}

Fix $u\in S^{d-1}$ and a nonempty $S\subseteq[2m]$ such that
$u\in\bigcap_{i\in S}G_i$. The boundary condition gives
\[
 \langle v_i,u\rangle=m_u:=\min_{j\in[2m]}\langle v_j,u\rangle
 \qquad(i\in S),
\]
where $m_u<0$ because $0\in\operatorname{int}P$. Thus the labels
$v_i$, $i\in S$, lie in the proper exposed face
$Q_u=\{z\in P:\langle z,u\rangle=m_u\}$, whose supporting hyperplane
avoids the origin. Since $P$ is simplicial, these labels are linearly independent. Since $v_{T(i)}=\pm v_i$, the relabeled vectors are also linearly
independent, so $0\notin\conv\{v_{T(i)}:i\in S\}$ and $T(S)$ is not
$(V,0)$-balanced.
As $T^{-1}=T$, this shows that the relabeled boundary cover
$(G_{T(i)}\cap S^{d-1})_{i=1}^{2m}$ has no $(V,0)$-balanced intersections.

Sending each $v_i$ to $v_{T(i)}$ and extending affinely over each face
defines a continuous odd map
$L_T:\partial P\to\mathbb{R}^d\setminus\{0\}$.
Define the radial projection
\[
 \tau:\mathbb{R}^d\setminus\{0\}\to S^{d-1},\qquad
 \tau(z)=\frac{z}{\|z\|}.
\]
Its restriction $\tau_P=\tau|_{\partial P}$ is an odd homeomorphism.
Consequently,
\[
 f_T=\tau\circ L_T\circ\tau_P^{-1}:S^{d-1}\to S^{d-1}
\]
is odd and hence is not null-homotopic by the Borsuk--Ulam theorem.

As in Definition~\ref{def:CoverDegree}, choose open neighborhoods
$G_i^\varepsilon\supseteq G_i\cap S^{d-1}$ in $S^{d-1}$ such that
$u\in G_i^\varepsilon$ implies $\langle u,v_i\rangle<0$ for
$i\in[2m]$, and, for every nonempty $S\subseteq[2m]$,
\[
 S^{d-1}\cap\bigcap_{i\in S}G_i\neq\varnothing
 \quad\Longleftrightarrow\quad
 \bigcap_{i\in S}G_i^\varepsilon\neq\varnothing.
\]
Let $(\psi_i)$ be a partition of unity subordinate to
$(G_i^\varepsilon)_{i=1}^{2m}$. The intersection condition ensures
that the labels with positive coefficients at each point lie in a
common proper face of $P$. Hence the map
$\rho(u)=\sum_{i=1}^{2m}\psi_i(u)v_i$ takes values in $\partial P$.
Since $\langle\rho(u),u\rangle<0$, the normalized straight-line
homotopy shows that $\tau_P\circ\rho$ is homotopic to the antipodal
map $u\mapsto-u$.

The functions $(\psi_{T(i)})$ form a partition of unity for the
relabeled boundary cover. Its normalized label map is
\[
 \begin{aligned}
 q_T(u)
 &=\tau\left(\sum_{i=1}^{2m}\psi_{T(i)}(u)v_i\right)\\
 &=\tau\bigl(L_T(\rho(u))\bigr)
 =f_T\bigl(\tau_P(\rho(u))\bigr).
 \end{aligned}
\]
Thus $q_T$ is homotopic to the odd map $u\mapsto f_T(-u)$ and has odd
degree. The relabeled boundary cover therefore has nonzero degree
relative to $(V,0)$. Theorem~\ref{CoverExtensionManifolds} 
forces a $(V,0)$-balanced intersection. It lies in the open
ball because the boundary cover has no such intersections.
\end{proof}

\section{Proof of Theorem~\ref{IntroLinDual}}

We use the notation of Theorem~\ref{IntroLinDual}. Weak duality follows
directly: if $Ax\geq b$, $\lambda\geq0$, and $A^T\lambda=c$, then
\[
 \langle x,c\rangle =\langle x,A^T\lambda\rangle =\langle\lambda,Ax\rangle \geq\langle\lambda,b\rangle.
\]

Before proving strong duality in the general case, we consider the
following special case.

\begin{theorem}
\label{TopLinDual}
Suppose that $0\notin\conv\{a_1,\ldots,a_m\}$, the rows $a_i$ span
$\mathbb{R}^n$, and $c\in\operatorname{cone}\{a_1,\ldots,a_m\}$.
Then the primal minimum and the dual maximum exist and are equal.
\end{theorem}

\begin{proof}
The case $n=1$ is immediate, so assume $n\geq2$.
Strict separation gives a unit vector $h$
with $\langle a_i,h\rangle>0$ for every $i$.
A sufficiently large multiple of $h$ is feasible; if $c=0$, take $\lambda=0$.
For $c\neq0$, positive rescaling of the inequalities and objective
allows us to assume
\[
 \langle a_i,h\rangle=\langle c,h\rangle=1\quad(i\in[m]), \qquad c\in\conv\{a_1,\ldots,a_m\}.
\]
Put $H=h^\perp$, $a_i'=a_i-h$, $c'=c-h$, and
$P=\conv\{a_1',\ldots,a_m'\}\subseteq H$.
Then $\dim P=n-1$ and $c'\in P$. For $x\in H$, set
\[
 \tau(x)=\max_i\bigl(b_i-\langle x,a_i'\rangle\bigr),\qquad
 F_i=\{x\in H:b_i-\langle x,a_i'\rangle=\tau(x)\}.
\]
These sets form a closed cover of $H$, and $z(x)=x+\tau(x)h$ is
feasible, with its $i$th constraint active exactly when $x\in F_i$.
Let $B_H$ be the closed unit ball of $H$, $S_H=\partial B_H$, and put
\[
 f(x)=\frac{x}{1+\|x\|},\qquad F_i'=\overline{f(F_i)}.
\]
The sets $F_i'$ cover $B_H$.

\begin{claim}
\label{BoundaryRule}
If $u\in F_i'\cap S_H$, then
\[
 \langle u,a_i'\rangle=\min_{j\in[m]}\langle u,a_j'\rangle.
\]
\end{claim}
\begin{proof}
Choose $x_k\in F_i$ with $\|x_k\|\to\infty$ and
$x_k/\|x_k\|\to u$. Maximality gives
$\langle a_i'-a_j',x_k\rangle\leq b_i-b_j$ for every $j\in[m]$.
Divide by $\|x_k\|$ and pass to the limit.
\end{proof}

For $c'$ interior to $P$ in $H$, apply Theorem~\ref{DualPolKKM} to
the cover $F_i'$, labels $a_i'$ and target $c'$.
By Claim~\ref{BoundaryRule}, a boundary intersection would have all
corresponding vectors $a_i'$ in a proper face of $P$, which cannot contain $c'$.
Thus the intersection lies in the open ball, where
$F_i'\cap(B_H\setminus S_H)=f(F_i)$, and yields a nonempty
$S\subseteq[m]$ such that
\[
 c'\in\conv\{a_i':i\in S\},\qquad
 \bigcap_{i\in S}F_i\neq\emptyset.
\]
For a boundary target, take interior points $c_k'\to c'$.
The preceding conclusion gives sets $S_k$; passing to a subsequence
makes $S_k=S$ fixed. Closedness of $\conv\{a_i':i\in S\}$ then gives
the same conclusion for $c'$.

Choose $x^*\in\bigcap_{i\in S}F_i$ and convex weights $\lambda_i\geq0$
supported on $S$ and representing $c'$. Then $z^*=x^*+\tau(x^*)h$
satisfies
\[
 Az^*\geq b,\qquad
 A^T\lambda=c,\qquad
 \langle z^*,c\rangle=\sum_i\lambda_i b_i.
\]
Weak duality and undoing the positive rescalings prove the assertion.
\end{proof}

\begin{proof}[Proof of Theorem~\ref{IntroLinDual}]
Weak duality was proved above. If neither constraint system is feasible,
neither extremum exists. The remaining possibilities are as follows.

\medskip
\noindent\textit{Case 1: Only the primal constraints are feasible.}
Choose $x_0$ with $Ax_0\geq b$.
Since $c$ lies outside the closed cone generated by the rows $a_i$,
separation gives $d$ with
\[
 \langle a_i,d\rangle\geq0\quad(i\in[m]),\qquad
 \langle c,d\rangle<0.
\]
The feasible ray $x_0+Td$, $T\geq0$, makes the primal objective tend
to $-\infty$.

\medskip
\noindent\textit{Case 2: Only the dual constraints are feasible.}
Choose $\lambda^0\geq0$ with $A^T\lambda^0=c$ and consider
\[
 \min\{t:\langle a_i,x\rangle+t\geq b_i\ (i\in[m]),\ t\geq0\}.
\]
Its normals are $(a_i,1)$ and $(0,1)$, and its objective vector is
$(0,1)$. Thus Theorem~\ref{TopLinDual} applies in the span of these
normals. It gives an attained minimum $t_1>0$, since $t_1=0$ would
make the original primal feasible, and dual coefficients satisfying
\[
 \eta\geq0,\qquad
 A^T\eta=0,\qquad
 \langle\eta,b\rangle=t_1>0.
\]
Hence the feasible dual ray $\lambda^0+T\eta$, $T\geq0$, makes the
dual objective tend to $+\infty$.

\medskip
\noindent\textit{Case 3: Both constraint systems are feasible.}
Choose $x_0$ and $\lambda^0\geq0$ with $Ax_0\geq b$ and
$A^T\lambda^0=c$. Choose $M>\sum_i\lambda_i^0$ so that $c/M$ lies
in no convex hull of rows that excludes zero. This is possible because
there are finitely many such convex hulls, each at positive distance
from zero.
Using the same augmented normals as in Case~2, apply
Theorem~\ref{TopLinDual} to
\[
 \min\{\langle x,c\rangle+Mt:
       \langle a_i,x\rangle+t\geq b_i\ (i\in[m]),\ t\geq0\}.
\]
The objective vector $(c,M)$ belongs to their cone, with coefficients
$(\lambda^0,M-\sum_i\lambda_i^0)$.
The resulting optimizers $(x^*,t^*)$ and $(\lambda,\nu)$, where $\nu$
corresponds to $t\geq0$, satisfy
\[
 A^T\lambda=c,\qquad \sum_i\lambda_i+\nu=M,\qquad
 \langle x^*,c\rangle+Mt^*=\langle\lambda,b\rangle.
\]
Equality in weak duality gives complementary slackness. If $t^*>0$,
then $\nu=0$, and $S=\{i:\lambda_i>0\}$ satisfies
\[
 c/M\in\conv\{a_i:i\in S\},\qquad
 \langle a_i,x^*\rangle+t^*=b_i\quad(i\in S).
\]
Our choice of $M$ implies $0\in\conv\{a_i:i\in S\}$, whereas
$\langle a_i,x_0-x^*\rangle\geq t^*>0$ for every $i\in S$, a
contradiction. Hence $t^*=0$, and $(x^*,\lambda)$ is feasible for the
original problems with equal objective values. Weak duality completes
the proof.
\end{proof}

\section{Proof of Theorem~\ref{IntroNonlinear}}
\label{sec:ProofNonlinear}

\begin{proof}[Proof of Theorem~\ref{IntroNonlinear}]
\leavevmode\par
\medskip
\noindent\textbf{Part 1: Feasibility and balancing coefficients.}\par\nobreak
\smallskip
If $c=0$, take all coefficients to be zero and choose a sufficiently
large multiple of $h$ as a feasible point. This proves the first
assertion in this case. For the rest of Part~1, assume $c\neq0$.
After an orthogonal change of coordinates, we may assume
$h=(0,\ldots,0,1)$.

Put $\alpha_i=\langle a_i,h\rangle>0$ for $i\in[m]$ and
$\beta=\langle c,h\rangle>0$, where the latter inequality follows from
$c\in\operatorname{cone}(A)\setminus\{0\}$.
Define
\[
 \widetilde a_i=\frac{a_i}{\alpha_i},\qquad
 \widetilde b_i=\frac{b_i}{\alpha_i},\qquad
 \widetilde r_i=\frac{r_i}{\alpha_i},\qquad
 \widetilde c=\frac{c}{\beta}.
\]
Dividing the entire $i$th constraint by $\alpha_i$ gives the equivalence
\[
 \langle a_i,x\rangle\geq b_i+r_i(\pi_Hx)
 \quad\Longleftrightarrow\quad
 \langle\widetilde a_i,x\rangle
 \geq\widetilde b_i+\widetilde r_i(\pi_Hx).
\]
Thus the feasible set and active indices are preserved.

Since all scaling factors are positive, for every $S\subseteq[m]$,
\[
 \begin{aligned}
 \operatorname{cone}\{a_i:i\in S\}
 &=\operatorname{cone}\{\widetilde a_i:i\in S\},\\
 c\in\operatorname{cone}\{a_i:i\in S\}
 &\quad\Longleftrightarrow\quad
 \widetilde c\in\operatorname{cone}\{\widetilde a_i:i\in S\}.
 \end{aligned}
\]
In particular, the generating property of $J$ is preserved. For each
$i\in[m]$, we also have
\[
 \begin{aligned}
 r_i(y)=o(\|y\|)
 &\quad\Longleftrightarrow\quad \widetilde r_i(y)=o(\|y\|),\\
 \max\{r_i(y),0\}=o(\|y\|)
 &\quad\Longleftrightarrow\quad
 \max\{\widetilde r_i(y),0\}=o(\|y\|).
 \end{aligned}
\]
Hence the asymptotic conditions are preserved as well.

Dropping the tildes for the rescaled data, we may write
\[
 a_i=(\bar a_i,1),\qquad
 c=(\bar c,1),\qquad
 x=(y,t).
\]
Identify $H$ with $\mathbb{R}^{n-1}$, so that
$\pi_Hx=y$. Put $\bar A=(\bar a_1,\ldots,\bar a_m)$ and
\[
 P=\conv(\bar A),\qquad \bar c \in P,\qquad
 \phi_i(y)=b_i+r_i(y)-\langle \bar a_i,y\rangle.
\]
As in the linear proof, $\dim P=n-1$.
Since all normalized rows have last coordinate $1$, the generating
property of $J$ also gives
\[
 P=\conv\{\bar a_j:j\in J\}.
\]

In these coordinates, the perturbed feasibility inequalities
$\langle x,a_i\rangle\geq b_i+r_i(\pi_Hx)$ become
$\langle(y,t),(\bar a_i,1)\rangle\geq b_i+r_i(y)$,
or equivalently, $t\geq\phi_i(y)$ for every $i\in[m]$.

We seek $y^*\in\mathbb{R}^{n-1}$ and coefficients $\lambda_i\geq0$ satisfying
\begin{equation}\label{GraphBalance}
 \sum_i\lambda_i=1,\qquad \sum_i\lambda_i \bar a_i=\bar c,\qquad
 \lambda_i>0\Longrightarrow
 \phi_i(y^*)=\max_j\phi_j(y^*).
\end{equation}
If $n=1$, take $y^*=0$ and assign coefficient $1$ to a maximizing index,
with all other coefficients zero. Now assume that $n\geq2$.

\smallskip
\noindent\textit{Case 1: $\bar c$ is interior to $P$.}
The sets
\[
 F_i=\{y:\phi_i(y)=\max_j\phi_j(y)\}
\]
form a closed cover of $\mathbb{R}^{n-1}$.
Let $B_H$ be the closed unit ball of $H$. As in the linear proof, put
\[
 f(y)=\frac{y}{1+\|y\|},\qquad F_i'=\overline{f(F_i)}\subseteq B_H.
\]
These sets form a closed cover of $B_H$.
If $u\in F_i'\cap\partial B_H$, choose $y_k\in F_i$ with
$\|y_k\|\to\infty$ and $y_k/\|y_k\|\to u$.
For every $j\in J$, maximality gives
\[
 \langle\bar a_i-\bar a_j,y_k\rangle
 \leq b_i-b_j+r_i(y_k)-r_j(y_k).
\]
Divide by $\|y_k\|$ and use
$\max\{r_i(y_k),0\}=o(\|y_k\|)$ and $r_j(y_k)=o(\|y_k\|)$.
We obtain
\[
 \langle u,\bar a_i\rangle\leq\langle u,\bar a_j\rangle \qquad(j\in J).
\]
Since $\bar a_i\in P=\conv\{\bar a_j:j\in J\}$, it follows that
\[
 \langle u,\bar a_i\rangle
 =\min_{j\in J}\langle u,\bar a_j\rangle
 =\min_{1\leq\ell\leq m}\langle u,\bar a_\ell\rangle.
\]
Thus the boundary rule of Claim~\ref{BoundaryRule} still holds.
Apply Theorem~\ref{DualPolKKM} to the cover $F_i'$ with vectors
$\bar a_i$ and target $\bar c$. This gives a nonempty $S\subseteq[m]$
and a point $z^*$ such that
\[
 z^*\in\bigcap_{i\in S}F_i',\qquad
 \bar c\in\conv\{\bar a_i:i\in S\}.
\]
As in the linear proof, $z^*$ lies in the open ball because $\bar c$
is interior to $P$. Define
\[
 y^*=\frac{z^*}{1-\|z^*\|}.
\]
Then $f(y^*)=z^*$. Since $F_i$ is closed,
$F_i'\cap\operatorname{int}B_H=f(F_i)$, so
$y^*\in\bigcap_{i\in S}F_i$.
Thus $y^*$ is an $(\bar A,\bar c)$-balanced intersection point of the original cover.
This gives \eqref{GraphBalance}.

\smallskip
\noindent\textit{Case 2: $\bar c$ is on the boundary of $P$.}
Under the normalization above, this is equivalent to
$c\in\partial\operatorname{cone}(A)$.
Choose $\bar c_k\in\operatorname{int}P$ with $\bar c_k\to\bar c$.
The sets $F_i$ are unchanged and do not depend on $c$. By Case~1, for each $k$
there is a nonempty $S_k\subseteq[m]$ such that
\[
 \bar c_k\in\conv\{\bar a_i:i\in S_k\},\qquad
 \bigcap_{i\in S_k}F_i\neq\emptyset.
\]
There are only finitely many index sets, so, after passing to a
subsequence, we may assume that $S_k=S$ is fixed.
Since $\conv\{\bar a_i:i\in S\}$ is closed, it contains $\bar c$.
Choose $y^*\in\bigcap_{i\in S}F_i$ and express $\bar c$ as a convex
combination of the vectors $\bar a_i$, $i\in S$.
Setting all remaining coefficients to zero gives \eqref{GraphBalance}.

By Carath\'eodory's theorem, the number of positive $\lambda_i$
can be reduced to at most $n$. Set
\[
 x^*=\bigl(y^*,\max_i\phi_i(y^*)\bigr).
\]
Then $\sum_i\lambda_i a_i=(\bar c,1)=c$.
By the choice of height, every perturbed inequality holds, with
equality whenever $\lambda_i>0$. For the original data, replace
$\lambda_i$ by $\beta\lambda_i/\alpha_i$, preserving the support of the coefficients,
and undo the orthogonal change of coordinates. This proves the first
assertion for the original data.

\medskip
\noindent\textbf{Part 2: Lagrangian lower bound.}\par\nobreak
\smallskip
We now use the original, unnormalized data, allowing $c=0$ as well.
Fix coefficients $\lambda_i\geq0$ with
$\sum_i\lambda_i a_i=c$. For every
$x\notin\bigcup_i\operatorname{int}(U_i)$,
\[
 \begin{aligned}
 &\langle x,c\rangle-r(\pi_Hx)\\
 &\quad\geq\sum_i\lambda_i b_i+
                 \sum_i\lambda_i r_i(\pi_Hx)-r(\pi_Hx)\\
 &\quad\geq\sum_i\lambda_i b_i+
       \inf_{y\in H}\left(\sum_i\lambda_i r_i(y)-r(y)\right).
 \end{aligned}
\]
Taking the infimum over all such $x$ and then the supremum over
all such coefficients proves the assertion. The argument also applies
when an inner infimum is $-\infty$.

\medskip
\noindent\textbf{Part 3: Scarf-type feasibility.}\par\nobreak
\smallskip
Take a feasible point $x^*$ and balancing coefficients for $(A,c)$ from
the first assertion, with at most $n$ positive entries. Complementary slackness
gives $x^*\in\bigcap_{\lambda_i>0}U_i$.
The assumed inclusion therefore implies $x^*\in D_b$.
This proves the feasibility claim and shows why it suffices to check
supports of size at most $n$.
\end{proof}

\begin{remark}
Taking $r_i=0$ in the first assertion of Theorem~\ref{IntroNonlinear}
gives Theorem~\ref{TopLinDual} by weak linear duality.
If $r_i\leq0$, then
$U_i\subseteq\{x:\langle a_i,x\rangle\leq b_i\}$. 
\end{remark}

\section{Proofs of Theorems~\ref{ManifoldComplementarity}
and~\ref{SymmetricPerturbations}}
\label{sec:ProofManifold}

\begin{proof}[Proof of Theorem~\ref{ManifoldComplementarity}]
Put $V=(v_1,\ldots,v_m)$.
The sets $F_i=\{x\in M:i\in I(x)\}$ form a closed cover of $M$.
The boundary condition excludes $(V,0)$-balanced intersections on $N$.
Moreover, nonzero degree makes $g$ surjective, so for every
$u\in S^{n-1}$ some label satisfies $\langle v_i,u\rangle<0$.
Hence $0\in\operatorname{int}\conv\{v_1,\ldots,v_m\}$.

Choose a partition of unity $(\psi_i)$ subordinate to a sufficiently
small open enlargement of $(F_i\cap N)$, preserving its intersecting
index sets and the strict boundary inequalities. Its normalized label map
\[
 q(x)=\frac{\sum_i\psi_i(x)v_i}{\|\sum_i\psi_i(x)v_i\|}
 \qquad(x\in N)
\]
satisfies $\langle q(x),g(x)\rangle<0$, so the normalized straight-line
homotopy joins $q$ to $-g$. Since the antipodal map has nonzero degree,
\[
 \deg_{V,0}\bigl((F_i\cap N)_{i=1}^m\bigr)=\deg(-g)\neq0.
\]
Theorem~\ref{CoverExtensionManifolds} therefore forces a $(V,0)$-balanced
intersection point $x^*\in M$ of $(F_i)$.
The strict boundary inequalities imply $x^*\notin N$.
Choose convex coefficients $\lambda_i$ representing $0$ among the labels indexed by
$I(x^*)$, and set all remaining coefficients to zero. Then
$\lambda_i\geq0$ for $i\in[m]$, $\sum_{i=1}^m\lambda_i=1$, and
$\sum_{i=1}^m\lambda_i v_i=0$.
Set $t^*=\Phi(x^*)$. Then $(x^*,t^*)$ is feasible, and
$t^*=\phi_i(x^*)$ whenever $\lambda_i>0$, proving complementary slackness.
\end{proof}

\medskip

\begin{remark}

The first assertion of Theorem~\ref{IntroNonlinear} can also be deduced
from Theorem~\ref{ManifoldComplementarity}. In the normalized setting,
with $n\geq2$, take $M=B_R\subseteq\mathbb{R}^{n-1}$ and restrict the
previously defined active cover to $M$. For interior targets, the
assumptions on $r_i$ ensure that the boundary condition holds with
$g_R(y)=y/R$ for every sufficiently large $R$.
With labels $v_i=\bar a_i-\bar c$ and target $0$, the boundary degree
therefore stabilizes at the nonzero degree of the antipodal map,
and Theorem~\ref{ManifoldComplementarity}
gives the required feasible point and balancing coefficients for $(A,c)$.
Boundary targets follow by approximation.
\end{remark}

\begin{proof}[Proof of Theorem~\ref{SymmetricPerturbations}]
Regroup the graph terms according to their linear parts by setting
\[
 \psi_i(y)=\max\{f_i^+(y),f_{i+m}^-(y)\}-\langle p_i,y\rangle
 \qquad(i\in[2m]).
\]

Regrouping preserves the upper envelope, so
$\max_{j\in[2m]}\psi_j(y)=\max_{j\in[2m]}\phi_j(y)=\Phi(y)$.

Define the permutation $T$ by
\[
 \bigl(T(i),T(i+m)\bigr)=
 \begin{cases}
 (i,i+m),
 &\text{if }f_i^-\leq f_{i+m}^+
       \text{ and }f_{i+m}^-\leq f_i^+,\\
 (i+m,i),&\text{otherwise}.
 \end{cases}
 \qquad(i\in[m]).
\]
The dominance condition gives $\psi_{T(i)}(y)\leq\phi_i(y)$ for every
$i\in[2m]$ and $y\in H$.

The sets $G_i=\{y\in H:\psi_i(y)=\Phi(y)\}$ form a closed cover of $H$.
Let $B_H$ be the closed unit ball of $H$ and put
\[
 \widehat G_i=\overline{\{y/(1+\|y\|):y\in G_i\}}\subseteq B_H
 \qquad(i\in[2m]).
\]
These sets form a closed cover of $B_H$. Since
$\psi_i(y)=-\langle p_i,y\rangle+o(\|y\|)$, the limiting argument
from the proof of Theorem~\ref{IntroNonlinear} gives
\[
 u\in\widehat G_i\cap\partial B_H
 \quad\Longrightarrow\quad
 \langle p_i,u\rangle=\min_{j\in[2m]}\langle p_j,u\rangle.
\]

We apply Theorem~\ref{AntipodalCover} with $d=n-1$, identifying $H$
isometrically with $\mathbb{R}^{n-1}$. Its ball is $B_H$, its cover is
$(\widehat G_i)_{i=1}^{2m}$, its label family is $V=(p_i)_{i=1}^{2m}$,
and its permutation is $T$. The label polytope $P$ is centrally
symmetric and simplicial, with the same antipodal pairs
$\{i,i+m\}$, $i\in[m]$.
The theorem therefore gives a nonempty $(V,0)$-balanced set $S\subseteq[2m]$
and a point $z^*\in\operatorname{int}(B_H)$ such that
\[
 z^*\in\bigcap_{i\in S}\widehat G_{T(i)}.
\]
Since each $G_i$ is closed in $H$, the inverse radial compactification
gives
\[
 y^*=\frac{z^*}{1-\|z^*\|}\in\bigcap_{i\in S}G_{T(i)}.
\]
Choose coefficients $\mu_i\geq0$, $i\in[2m]$, with $\mu_i=0$ for
$i\notin S$, $\sum_{i=1}^{2m}\mu_i=1$, and
$\sum_{i=1}^{2m}\mu_i p_i=0$.
Then $\psi_{T(i)}(y^*)=\Phi(y^*)$ whenever $\mu_i>0$.
By Carath\'eodory's theorem, at most $n$ of the coefficients need be positive.

Set $t^*=\Phi(y^*)$ and $x^*=y^*+t^*h$. This point is feasible.
If $\mu_i>0$, then
\[
 \Phi(y^*)=\psi_{T(i)}(y^*)\leq\phi_i(y^*)\leq\Phi(y^*),
\]
so $t^*=\phi_i(y^*)$ and $x^*\in\partial U_i$.
Finally, the coefficients $\lambda_i=\mu_i/\alpha_i$ have the same
support and satisfy
\[
 \sum_{i=1}^{2m}\lambda_i a_i
 =\sum_{i=1}^{2m}\mu_i(h+p_i)=h=c.
\]
\end{proof}

\section{Cooperative games with non-transferable utility}
\label{sec:NTUGames}

In this section we derive Scarf's theorem directly from
Theorem~\ref{IntroNonlinear}. We begin with the definition of a cooperative game.

Let $N=\{1,\ldots,n\}$ be the set of players; its nonempty subsets
are called coalitions. For every
coalition $S\subseteq N$, identify $\mathbb{R}^S$ with the coordinate
subspace of $\mathbb{R}^n$ consisting of vectors whose coordinates
outside $S$ are zero. Write $\mathbf{1}_S$ for the characteristic
vector of $S$ and $x_S$ for the restriction of $x$ to the coordinates
in $S$.

\begin{definition}[Normalized cooperative NTU game]
\label{def:NormalizedNTU}
A normalized cooperative NTU game is a triple $(N,U,V)$ such that, for
each coalition $S$, the set $V(S)\subseteq\mathbb{R}^S$ is nonempty,
closed and \emph{comprehensive}: $v\in V(S)$ and $w\in\mathbb{R}^S$
with $w\leq v$ imply $w\in V(S)$.
The sets $U(S)=V(S)\times\mathbb{R}^{N\setminus S}$ are cylindrical.
There is a common $M>0$ such that
$\langle v,\mathbf{1}_S\rangle\leq M$ for all coalitions $S$ and $v\in V(S)$,
and normalization requires
\[
 V(\{i\})=\{v\in\mathbb{R}^{\{i\}}:v_i\leq0\}\qquad(i\in N).
\]
\end{definition}

An allocation is a vector $x\in V(N)$. A coalition $S$ \emph{blocks}
$x$ if there exists $v\in V(S)$ such that $v_S>x_S$.

\begin{definition}
The \emph{core} is the set of allocations that are not blocked by any
coalition:
\[
 C(U)=\{x\in V(N):\text{no nonempty coalition }S\subseteq N
                         \text{ blocks }x\}.
\]
\end{definition}

The game is \emph{balanced} if, for every family of coefficients
$\lambda_S\geq0$ satisfying
\[
 \sum_{\emptyset\neq S\subseteq N}\lambda_S\mathbf{1}_S=\mathbf{1}_N,
\]
we have
\[
 \bigcap_{\lambda_S>0}U(S)\subseteq U(N).
\]

Thus the condition says that feasibility for every coalition with
positive weight implies feasibility for the grand coalition.
We now derive Scarf's core existence theorem directly from the third
assertion of Theorem~\ref{IntroNonlinear}.

\begin{remark}
Our boundedness assumption is stronger than the standard one, which
requires only that $V(S)\cap\mathbb{R}_{\geq0}^S$ be bounded for every
coalition $S$ \cite[Definition~1.5]{BludovMusin}.
Since individual payoffs are normalized to zero, any allocation with a
negative coordinate is blocked by a singleton. Moreover, a coalition
blocking a nonnegative allocation must offer a payoff vector that is
positive in every coordinate of that coalition. Thus the core depends
only on the nonnegative parts $V(S)\cap\mathbb{R}_{\geq0}^S$; payoff
vectors outside these parts are irrelevant to core membership. Thus
Scarf's theorem for our restricted class of games immediately implies
Scarf's theorem for the standard class of games.

\end{remark}

\begin{theorem}[Scarf~\cite{Scarf}]
\label{ScarfNTU}
Every balanced normalized cooperative NTU game has a nonempty core.
\end{theorem}

\begin{proof}

Put $h=\mathbf{1}_N/\sqrt n$ and $H=h^\perp$. For every nonempty
coalition $S\subseteq N$ and $y\in H$, define
\[
 \tau_S(y)=\max\{t\in\mathbb{R}:y+t\mathbf{1}_N\in U(S)\},\qquad
 r_S(y)=\langle\mathbf{1}_S,y\rangle+|S|\tau_S(y).
\]
The assumptions on $V(S)$ ensure that the admissible heights form a
nonempty closed lower ray bounded above, so $\tau_S(y)$ is finite.
Comprehensiveness gives
$|\tau_S(y)-\tau_S(z)|\leq\|y-z\|_\infty$, hence $r_S$ is continuous.
Moreover,
\[
 U(S)=\{x:\langle\mathbf{1}_S,x\rangle\leq r_S(\pi_Hx)\}, \qquad r_S(y)\leq M.
\]
Thus $\max\{r_S(y),0\}=o(\|y\|)$, while normalization gives
$r_{\{i\}}=0$. The hypotheses of
Theorem~\ref{IntroNonlinear} hold for $a_S=\mathbf{1}_S$ and $b_S=0$,
with $J$ consisting of the singleton coalitions, whose vectors generate
the cone $\mathbb{R}_{\geq0}^n$ of all coalition vectors.

Take $c=\mathbf{1}_N$, $r=r_N$ and $b=0$, so that $D_0=U(N)$.
Balancedness is precisely the inclusion required in Part~3 of
Theorem~\ref{IntroNonlinear}, which therefore gives
\[
 x^*\in U(N)\setminus
       \bigcup_{\emptyset\neq S\subseteq N}\operatorname{int}U(S).
\]
If a coalition $S$ blocked $x^*$, comprehensiveness and cylindricity
would imply $x^*\in\operatorname{int}U(S)$, a contradiction.
Thus $x^*\in C(U)$.
\end{proof}

\section{Discussion and conclusion}
\label{sec:Conclusion}

We propose a general KKM-type covering framework for finding feasible
points and balancing coefficients supported on active constraints.
It gives a direct proof of
linear duality and extends it to continuously perturbed constraints,
yielding Scarf's theorem for normalized cooperative games with
non-transferable utility. Theorem~\ref{ManifoldComplementarity} gives
complementary slackness on compact connected oriented manifolds under
a nonzero-degree boundary condition. This suggests a way to model
further problems with manifold-valued parameters and balance
requirements on the active constraints.

Theorem~\ref{SymmetricPerturbations} relies on antipodal symmetry and
the Borsuk--Ulam theorem. Following the equivariant approach of
Avvakumov and Karasev~\cite{Avvakumov-Karasev} to envy-free division
with possibly degenerate pieces, $\mathbb{Z}_{p^k}$-equivariant analogues
of Theorem~\ref{AntipodalCover} and, consequently,
Theorem~\ref{SymmetricPerturbations} follow by the same argument under
appropriate symmetry and dominance assumptions, where $p$ is prime
and $k\geq1$.
The corresponding equivariant degree result replaces the
Borsuk--Ulam theorem.

\section*{Acknowledgments}

The author is grateful to Dmitry Gribanov, Dmitry Kovalenko, Alexander Levin,
and the International Laboratory of Game Theory and Decision Making for helpful discussions. The author is grateful to Georgy Potapov and Roman Karasev for reading the initial draft of the manuscript. The work presented in Section~3 was carried out within the framework of the Basic Research Program at HSE University (HSE-BR-2025-077).

\section*{Declaration on the use of AI}

The main mathematical ideas in this paper were developed by the author.
GPT 6.1 was used to assist with developing, drafting, and checking parts
of proofs, searching the literature, and editing the language. 
The final version of Theorem~\ref{IntroNonlinear} was found during the discussion with AI.
The author is responsible for the final content of the paper.

 \medskip
 
Mikhail Bludov, Moscow Institute of Physics and Technology, Dolgoprudny, Moscow
region, Russia
 
  {\it E-mail address:} bludov.mv@phystech.edu

\end{document}